\def\coralreport{1}

\documentclass[10pt]{article}

\usepackage{ifthen}
\ifthenelse{\coralreport = 1}{
\usepackage{isetechreport}
\usepackage{amsthm}
\newtheorem{theorem}{Theorem}[section]
\newtheorem{lemma}[theorem]{Lemma}

\def\reportyear{16T}
\def\reportno{010}
\def\revisionno{0}
\def\originaldate{\today}
\def\revisiondate{???}
\coraltrue
\cvcrfalse
\usepackage{hyperref}
}{}

\usepackage{fec}
\usepackage{graphicx}

\begin{document}

\title{$R$-Linear Convergence of Limited Memory Steepest Descent\thanks{This material is based upon work supported by the U.S.~Department of Energy, Office of Science, Office of Advanced Scientific Computing Research, Applied Mathematics, Early Career Research Program under contract number DE--SC0010615, as well as by the U.S.~National Science Foundation under Grant Nos.~DMS--1319356 and CCF-1618717.}}

\ifthenelse{\coralreport = 1}{
\author{Frank E. Curtis\thanks{E-mail: \texttt{\href{mailto:frank.e.curtis@gmail.com}{frank.e.curtis@gmail.com}}}}
\author{Wei Guo\thanks{E-mail: \texttt{\href{mailto:weg411@lehigh.edu}{weg411@lehigh.edu}}}}
\affil{Department of Industrial and Systems Engineering, Lehigh University, USA}
\titlepage
}{
\author       {{\sc Frank~E.~Curtis\thanks{E-mail: frank.e.curtis@gmail.com} and Wei Guo\thanks{E-mail: weg411@lehigh.edu}} \\[2pt]
                 Dept. of Industrial and Systems Engineering, Lehigh University, Bethlehem, PA, USA}
\shortauthorlist{F.~E.~Curtis and W.~Guo}
\date         {\today}
}

\maketitle

\begin{abstract}
  {The limited memory steepest descent method (LMSD) proposed by Fletcher is an extension of the Barzilai-Borwein ``two-point step size'' strategy for steepest descent methods for solving unconstrained optimization problems.  It is known that the Barzilai-Borwein strategy yields a method with an $R$-linear rate of convergence when it is employed to minimize a strongly convex quadratic.  This paper extends this analysis for LMSD, also for strongly convex quadratics.  In particular, it is shown that the method is $R$-linearly convergent for any choice of the history length parameter.  The results of numerical experiments are provided to illustrate behaviors of the method that are revealed through the theoretical analysis.}
  {
    \ifthenelse{\coralreport = 1}{\\ Keywords:}{}
    unconstrained optimization;
    steepest descent methods;
    Barzilai-Borwein methods;
    limited memory methods;
    quadratic optimization;
    $R$-linear rate of convergence
  }
\end{abstract}

\bcenter\textit{Dedicated to Roger Fletcher and Jonathan Borwein whose contributions\\ continue to inspire many in the fields of nonlinear optimization and applied mathematics}\ecenter

\section{Introduction}\label{sec.introduction}

For solving unconstrained nonlinear optimization problems, one of the simplest and most widely used techniques is \emph{steepest descent} (SD).  This refers to any strategy in which, from any solution estimate, a productive step is obtained by moving some distance along the negative gradient of the objective function, i.e., the direction along which function descent is steepest.

While SD methods have been studied for over a century and employed in numerical software for decades, a unique and powerful instance came about relatively recently in the work by \cite{BarzBorw88}, where a ``two-point step size'' strategy is proposed and analyzed.  The resulting SD method, commonly referred to as the BB method, represents an effective alternative to other SD methods that employ an exact or inexact line search when computing the stepsize in each iteration.

The theoretical properties of the BB method are now well-known when it is employed to minimize an $n$-dimensional strongly convex quadratic objective function.  Such objective functions are interesting in their own right, but one can argue that such analyses also characterize the behavior of the method in the neighborhood of a strong local minimizer of any smooth objective function.  In the original work (i.e., \cite{BarzBorw88}), it is shown that the method converges $R$-superlinearly when $n=2$.  In \cite{Rayd93}, it is shown that the method converges from any starting point for any natural number~$n$, and in \cite{DaiLiao02} it is shown that the method converges $R$-linearly for any such $n$.

In each iteration of the BB method, the stepsize is determined by a computation involving the displacement in the gradient of the objective observed between the current iterate and the previous iterate.  As shown in \cite{Flet12}, this idea can be extended to a \emph{limited memory steepest descent} (LMSD) method in which a \emph{sequence} of $m$ stepsizes is computed using the displacements in the gradient over the previous $m$ steps.  This extension can be motivated by the observation that these displacements lie in a Krylov subspace determined by a gradient previously computed in the algorithm, which in turn yields a computationally efficient strategy for computing $m$ distinct eigenvalue estimates of the Hessian (i.e., matrix of second derivatives) of the objective function.  The reciprocals of these eigenvalue estimates represent reasonable stepsize choices.  Indeed, if the eigenvalues of the Hessian are computed exactly, then the algorithm terminates in a finite number of iterations; see \cite{Flet12} and \S\ref{sec.fundamentals}.

In \cite{Flet12}, it is shown that the proposed LMSD method converges from any starting point when it is employed to minimize a strongly convex quadratic function.  However, to the best of our knowledge, the convergence rate of the method for $m > 1$ has not yet been analyzed.  The main purpose of this paper is to show that this LMSD method converges $R$-linearly when employed to minimize such a function.  Our analysis builds upon the analyses in \cite{Flet12} and \cite{DaiLiao02}.

We mention at the outset that numerical evidence has shown that the practical performance of the BB method is typically much better than known convergence proofs suggest; in particular, the empirical rate of convergence is often $Q$-linear with a contraction constant that is better than that observed for a basic SD method.  Based on such evidence, we do not claim that the convergence results proved in this paper fully capture the practical behavior of LMSD methods.  To explore this claim, we present the results of numerical experiments that illustrate our convergence theory and demonstrate that the practical performance of LMSD can be even better than the theory suggests.  We conclude with a discussion of possible explanations of why this is the case for LMSD, in particular by referencing a known finite termination result for a special (computationally expensive) variant of the algorithm.

\paragraph{Organization}

In \S\ref{sec.fundamentals}, we formally state the problem of interest, notation to be used throughout the paper, Fletcher's LMSD algorithm, and a finite termination property for it.  In~\S\ref{sec.linear}, we prove that the LMSD algorithm is $R$-linearly convergent for any history length.  The theoretical results proved in~\S\ref{sec.linear} are demonstrated numerically in \S\ref{sec.numerical} and concluding remarks are presented in~\S\ref{sec.conclusion}.

\paragraph{Notation}

The set of real numbers (i.e., scalars) is denoted as $\R{}$, the set of nonnegative real numbers is denoted as $\R{}_+$, the set of positive real numbers is denoted as $\R{}_{++}$, and the set of natural numbers is denoted as $\N{} := \{1,2,\dots\}$.  A natural number as a superscript is used to denote the vector-valued extension of any of these sets---e.g., the set of $n$-dimensional real vectors is denoted as $\R{n}$---and a Cartesian product of natural numbers as a superscript is used to denote the matrix-valued extension of any of these sets---e.g., the set of $n \times n$ real matrices is denoted as $\R{n\times n}$.  A finite sequence of consecutive positive integers of the form $\{1,\dots,n\} \subset \N{}$ is denoted using the shorthand $[n]$.  Subscripts are used to refer to a specific element of a sequence of quantities, either fixed or generated by an algorithm.  For any vector $v \in \R{n}$, its Euclidean (i.e., $\ell_2$) norm is denoted by $\|v\|$.

\section{Fundamentals}\label{sec.fundamentals}

In this section, we state the optimization problem of interest along with corresponding definitions and concepts to which we will refer throughout the remainder of the paper.  We then state Fletcher's LMSD algorithm and prove a finite termination property for it similar to that proved in \cite{Flet12}.

\subsection{Problem Statement}

Consider the problem to minimize a strongly convex quadratic function $f : \R{n} \to \R{}$ defined by a symmetric positive definite matrix $A \in \R{n \times n}$ and vector $b \in \R{n}$, namely,
\bequation\label{prob.f}
  \min_{x\in\R{n}}\ f(x),\ \ \text{where}\ \ f(x) = \thalf x^TAx - b^Tx.
\eequation
Formally, we make the following assumption about the problem data.

\bassumption\label{ass.f}
  \textit{
  The matrix $A$ in problem~\eqref{prob.f} has $r \leq n$ distinct eigenvalues denoted by
  \bequation\label{eq.eigenvalues}
    \lambda_{(r)} > \cdots > \lambda_{(1)} > 0.
  \eequation
  Consequently, this matrix yields the eigendecomposition $A = Q \Lambda Q^T$, where
  \bequation\label{eq.eigendecomposition}
    \baligned
      Q &= \bbmatrix q_1 & \cdots & q_n \ebmatrix && \text{is orthogonal} \\
      \text{and}\ \ \Lambda &= \diag(\lambda_1,\dots,\lambda_n) && \text{with}\ \ \lambda_n > \cdots > \lambda_1 > 0 \\
      &&& \text{and}\ \ \lambda_i \in \{\lambda_{(1)},\dots,\lambda_{(r)}\}\ \ \text{for all}\ \ i \in [n].
    \ealigned
  \eequation
  }
\eassumption

The eigendecomposition of $A$ defined in Assumption~\ref{ass.f} plays a crucial role in our analysis.  In particular, we will make extensive use of the fact that any gradient of the objective function computed in the algorithm, a vector in $\R{n}$, can be written as a linear combination of the columns of the orthogonal matrix $Q$.  This will allow us to analyze the behavior of the algorithm componentwise according to the weights in these linear combinations corresponding to the sequence of computed objective gradients.  Such a strategy has been employed in all of the aforementioned articles on BB and LMSD.

\subsection{Limited memory steepest descent (LMSD) method}

Fletcher's limited memory steepest descent method is stated as Algorithm~\ref{alg.lmsd}.  The iterate update in the algorithm is the standard update in an SD method: each subsequent iterate is obtained from the current iterate minus a multiple of the gradient of the objective function evaluated at the current iterate.  With this update at its core, Algorithm~\ref{alg.lmsd} operates in cycles.  At $x_{k,1} \in \R{n}$ representing the initial point of the~$k$th cycle, a sequence of~$m$ positive stepsizes $\{\alpha_{k,j}\}_{j\in[m]}$ are selected to be employed in an inner cycle composed of $m$ updates, the result of which is set as the initial point for cycle $k+1$.

Once such an inner cycle has been performed, the stepsizes to be employed in the next cycle are computed as the reciprocals of estimates of eigenvalues of $A$.  \cite{Flet12} describes how these can be obtained in one of three ways, all offering the same estimates (in exact arithmetic).  The most intuitive definition is that, for cycle $k+1$, the estimates come as the eigenvalues of $T_k := Q_k^TAQ_k$, where $Q_k \in \R{n \times m}$ satisfying $Q_k^TQ_k = I$ is defined by a thin QR factorization of the matrix of $k$th cycle gradients, i.e., for some upper triangular matrix $R_k \in \R{m \times m}$, such a factorization satisfies the equation
\bequation\label{eq.G}
  Q_kR_k = G_k := \bbmatrix g_{k,1} & \cdots & g_{k,m} \ebmatrix.
\eequation
(For now, let us assume that $G_k$ has linearly independent columns, in which case the matrix $R_k$ in~\eqref{eq.G} is nonsingular.  For a discussion of situations when this is not the case, see Remark~\ref{rem.G_lin_ind} later on.)  Practically, however, obtaining $T_k$ in this manner requires multiplications with $A$ as well as the storage of the $n$-vectors composing the columns of $Q_k$.  Both can be avoided by obtaining the gradient at the initial point of cycle $k+1$, namely $g_{k+1,1} \equiv g_{k,m+1}$, as well as the matrix of $k$th-cycle reciprocal stepsizes
\bequation\label{eq.J}
  J_k \gets \bbmatrix \alpha_{k,1}^{-1} & & \\ -\alpha_{k,1}^{-1} & \ddots & \\ & \ddots & \alpha_{k,m}^{-1} \\ & & -\alpha_{k,m}^{-1} \ebmatrix,
\eequation
computing (upper triangular) $R_k$ and $r_k$ from the partially extended Cholesky factorization
\bequation\label{eq.Cholesky}
  G_k^T\bbmatrix G_k & g_{k,m+1} \ebmatrix = R_k^T\bbmatrix R_k & r_k \ebmatrix,
\eequation
then computing
\bequation\label{eq.T}
  T_k \gets \bbmatrix R_k & r_k \ebmatrix J_k R_k^{-1}.
\eequation
Fletcher's third approach, which also avoids multiplications with $A$, is to compute
\bequation\label{eq.QR}
  T_k \gets \bbmatrix R_k & Q_k^Tg_{k,m+1} \ebmatrix J_k R_k^{-1}.
\eequation
However, this is less efficient than using \eqref{eq.T} due to the need to store $Q_k$ and since the QR factorization of $G_k$ requires $\sim\!\!m^2 n$ flops, as opposed to the $\sim\!\!\thalf m^2 n$ flops required for \eqref{eq.T}; see \cite{Flet12}.

\balgorithm[H]
  \renewcommand{\thealgorithm}{LMSD}
  \caption{\ Limited Memory Steepest Descent Method}
  \label{alg.lmsd}
  \begin{algorithmic}[1]
    \State choose an initial point $x_{1,1} \in \R{n}$, history length $m \in \N{}$, and termination tolerance $\epsilon \in \R{}_{+}$
    \State choose stepsizes $\{\alpha_{1,j}\}_{j\in[m]} \subset \R{}_{++}$
    \State compute $g_{1,1} \gets \nabla f(x_{1,1})$
    \State \textbf{if} $\|g_{1,1}\| \leq \epsilon$, \textbf{then return} $x_{1,1}$ \label{step.termination1}
    \For{$k \in \N{}$}
      \For{$j \in [m]$}
        \State set $x_{k,j+1} \gets x_{k,j} - \alpha_{k,j} g_{k,j}$
        \State compute $g_{k,j+1} \gets \nabla f(x_{k,j+1})$
        \State \textbf{if} $\|g_{k,j+1}\| \leq \epsilon$, \textbf{then return} $x_{k,j+1}$ \label{step.termination2}
      \EndFor
      \State set $x_{k+1,1} \gets x_{k,m+1}$ and $g_{k+1,1} \gets g_{k,m+1}$ \label{step.end_of_cycle}
      \State set $G_k$ by \eqref{eq.G} and $J_k$ by \eqref{eq.J}
      \State compute $R_k$ and $r_k$ to satisfy \eqref{eq.Cholesky} and set $T_k$ by \eqref{eq.T}
      \State set $\{\theta_{k,j}\}_{j\in[m]} \subset \R{}_{++}$ as the eigenvalues of $T_k$ in decreasing order
      \State set $\{\alpha_{k+1,j}\}_{j\in[m]} \gets \{\theta_{k,j}^{-1}\}_{j\in[m]} \subset \R{}_{++}$
    \EndFor
  \end{algorithmic}
\end{algorithm}

The choice to order the eigenvalues of $T_k$ in decreasing order is motivated by \cite{Flet12}.  In short, this ensures that the stepsizes in cycle $k+1$ are ordered from smallest to largest, which improves the likelihood that the objective function and the norm of the objective gradient decrease monotonically, at least initially, in each cycle.  This ordering is not essential for our analysis, but is a good choice for any implementation of the algorithm; hence, we state the algorithm to employ this ordering.

One detail that remains for a practical implementation of the method is how to choose the initial stepsizes $\{\alpha_{1,j}\}_{j\in[m]} \subset \R{}_{++}$.  This choice has no effect on the theoretical results proved in this paper, though our analysis does confirm the fact that the practical performance of the method can improved if one has the knowledge to choose one or more stepsizes exactly equal to reciprocals of eigenvalues of~$A$; see~\S\ref{sec.finite}.  Otherwise, one can either provide a full set of $m$ stepsizes or carry out an initialization phase in which the first few cycles are shorter in length, dependent on the number of objective gradients that have been observed so far; see \cite{Flet12} for further discussion on this matter.

\bremark\label{rem.G_lin_ind}
  \textit{
  In \eqref{eq.G}, if $G_k$ for some $k \in \N{}$ does not have linearly independent columns, then $R_k$ is singular and the formulas \eqref{eq.T} and \eqref{eq.QR} are invalid, meaning that the employed approach is not able to provide~$m$ eigenvalue estimates for cycle~$k$.  As suggested in \cite{Flet12}, an implementation of the method can address this by iteratively removing ``older'' columns of $G_k$ until the columns form a linearly independent set of vectors, in which case the approach would be able to provide $\mtilde \leq m$ stepsizes for the subsequent (shortened) cycle.  We advocate such an approach in practice and, based on the results proved in this paper, conjecture that the convergence rate of the algorithm would be $R$-linear.  However, the analysis for such a method would be extremely cumbersome given that the number of iterations in each cycle might vary from one cycle to the next within a single run of the algorithm.  Hence, in our analysis in \S\ref{sec.linear}, we assume that $G_k$ has linearly independent columns for all $k \in \N{}$.  In fact, we go further and assume that $\|R_k^{-1}\|$ is bounded proportionally to the reciprocal of the norm of the objective gradient at the first iterate in cycle $k$ (meaning that the upper bound diverges as the algorithm converges to the minimizer of the objective function).  These norms are easily computed in an implementation of the algorithm; hence, we advocate that a procedure of iteratively removing ``older'' columns of $G_k$ would be based on observed violations of such a bound.  See the discussion following Assumption~\ref{ass.lmsd} in \S\ref{sec.linear}.
  }
\eremark

\subsection{Finite Termination Property of LMSD}\label{sec.finite}

If, for some $k \in \N{}$ and $j \in [m]$, the stepsizes in Algorithm~\ref{alg.lmsd} up through iteration $(k,j) \in \N{} \times [m]$ include the reciprocals of all of the $r \leq n$ distinct eigenvalues of~$A$, then the algorithm terminates by the end of iteration $(k,j)$ with $x_{k,j+1}$ yielding $\|g_{k,j+1}\| = 0$.  This is shown in the following lemma and theorem, which together demonstrate and extend the arguments made in \S2 of \cite{Flet12}.

\blemma\label{lem.recursion}
  \textit{
  Under Assumption~\ref{ass.f}, for each $(k,j) \in \N{} \times [m]$, there exist weights $\{d_{k,j,i}\}_{i\in[n]}$ such that~$g_{k,j}$ can be written as a linear combination of the columns of $Q$ in \eqref{eq.eigendecomposition}, i.e.,
  \bequation\label{eq.combination}
    g_{k,j} = \sum_{i=1}^n d_{k,j,i} q_i.
  \eequation
  Moreover, these weights satisfy the recursive property
  \bequation\label{eq.recursion}
    d_{k,j+1,i} = (1 - \alpha_{k,j} \lambda_i) d_{k,j,i}\ \ \text{for all}\ \ (k,j,i) \in \N{} \times [m] \times [n].
  \eequation
  }
\elemma 
\bproof
  Since $g_{k,j} = Ax_{k,j} - b$ for all $(k,j) \in \N{} \times [m]$, it follows that
  \bequalin
             &&  x_{k,j+1} &=  x_{k,j} - \alpha_{k,j}g_{k,j},\\
    \implies && Ax_{k,j+1} &= Ax_{k,j} - \alpha_{k,j}Ag_{k,j},\\
    \implies &&  g_{k,j+1} &=  g_{k,j} - \alpha_{k,j}Ag_{k,j},\\
    \implies &&  g_{k,j+1} &= (I - \alpha_{k,j}A)g_{k,j},\\
    \implies &&  g_{k,j+1} &= (I - \alpha_{k,j}Q \Lambda Q^T) g_{k,j},
  \eequalin
  from which one obtains that
  \bequationn
    \sum_{i=1}^n d_{k,j+1,i} q_i = \sum_{i=1}^n d_{k,j,i}(I - \alpha_{k,j} Q \Lambda Q^T)q_i = \sum_{i=1}^n d_{k,j,i}(q_i - \alpha_{k,j}\lambda_iq_i) = \sum_{i=1}^n d_{k,j,i}(1 - \alpha_{k,j}\lambda_i)q_i.
  \eequationn
  The result then follows since the columns of $Q$ form an orthogonal basis of $\R{n}$.
\end{proof}

\btheorem\label{th.finite}
  \textit{
  Suppose that Assumption~\ref{ass.f} holds and that Algorithm~\ref{alg.lmsd} is run with termination tolerance $\epsilon = 0$.  If, for some $(k,j) \in \N{} \times [m]$, the set of computed stepsizes up through iteration $(k,j)$ includes all of the values $\{\lambda_{(l)}^{-1}\}_{l\in[r]}$, then, at the latest, the algorithm terminates finitely at the end of iteration $(k,j)$ with $x_{k,j+1}$ yielding $\|g_{k,j+1}\| = 0$.
  }
\etheorem
\bproof
  Consider any $(k,j) \in \N{} \times [m]$ such that the stepsize is equal to the reciprocal of an eigenvalue of $A$, i.e., $\alpha_{k,j} = \lambda_{(l)}^{-1}$ for some $l \in [r]$.  By Lemma~\ref{lem.recursion}, it follows that
  \bequationn
    d_{k,j+1,i} = (1 - \alpha_{k,j} \lambda_i) d_{k,j,i} = (1 - \lambda_{(l)}^{-1} \lambda_i) d_{k,j,i} = 0\ \ \text{for all}\ \ i \in [n]\ \ \text{such that}\ \ \lambda_i = \lambda_{(l)}.
  \eequationn
  Along with the facts that Lemma~\ref{lem.recursion} also implies
  \bequationn
    d_{k,j,i} = 0 \implies d_{k,j+1,i} = 0\ \ \text{for all}\ \ (k,j) \in \N{} \times [m]
  \eequationn
  and $x_{k+1,1} \gets x_{k,m+1}$ (and $g_{k+1,1} \gets g_{k,m+1}$) for all $k \in \N{}$, the desired conclusion follows.
\eproof

\bremark
  \textit{
  Theorem~\ref{th.finite} implies that Algorithm~\ref{alg.lmsd} will converge finitely by the end of the second cycle if $m \geq r$ and the eigenvalues of $T_1$ include all eigenvalues $\{\lambda_{(l)}\}_{l\in[r]}$.  This is guaranteed, e.g., when the first cycle involves $m = n$ steps and $G_1$ has linearly independent columns.
  }
\eremark

\section{$R$-Linear Convergence Rate of LMSD}\label{sec.linear}

Our primary goal in this section is to prove that Algorithm~\ref{alg.lmsd} converges $R$-linearly for any choice of the history length parameter $m \in \N{}$.  For context, we begin by citing two known convergence results that apply for Algorithm~\ref{alg.lmsd}, then turn our attention to our new convergence rate results.

\subsection{Known Convergence Properties of LMSD}

In the Appendix of \cite{Flet12}, the following convergence result is proved for Algorithm~\ref{alg.lmsd}.  The theorem is stated slightly differently here only to account for our different notation.

\btheorem\label{th.Fletcher}
  \textit{
  Suppose that Assumption~\ref{ass.f} holds and that Algorithm~\ref{alg.lmsd} is run with termination tolerance $\epsilon = 0$.  Then, either $g_{k,j} = 0$ for some $(k,j) \in \N{} \times [m]$ or the sequences $\{g_{k,j}\}_{k=1}^\infty$ for each $j \in [m]$ converge to zero.
  }
\etheorem

\noindent
As a consequence of this result, we may conclude that if Algorithm~\ref{alg.lmsd} does not terminate finitely, then, according to the relationship \eqref{eq.combination}, the following limits hold:
\bsubequations\label{eq.limits}
  \begin{align}
    \lim_{k\to\infty}\ g_{k,j}   &= 0\ \ \text{for each}\ \ j \in [m]\ \ \text{and} \label{eq.limit_g} \\
    \lim_{k\to\infty}\ d_{k,j,i} &= 0\ \ \text{for each}\ \ (j,i) \in [m] \times [n]. \label{eq.limit_d}
  \end{align}
\esubequations
Fletcher's result, however, does not illuminate the rate at which these sequences converge to zero.  Only for the case of $m=1$ in which Algorithm~\ref{alg.lmsd} reduces to a BB method do the following results from \cite{DaiLiao02} (see Lemma~2.4 and Theorem~2.5 therein) provide a convergence rate guarantee.

\begin{lemma}\label{lem.DaiLiao}
  \textit{
  Suppose that Assumption~\ref{ass.f} holds and that Algorithm~\ref{alg.lmsd} is run with history length $m = 1$ and termination tolerance $\epsilon = 0$.  Then, there exists $K \in \N{}$, dependent only on $(\lambda_1,\lambda_n)$, such that
  \bequationn
    \|g_{k+K,1}\| \leq \thalf \|g_{k,1}\|\ \ \text{for all}\ \ k \in \N{}.
  \eequationn
  }
\end{lemma}

\btheorem\label{th.DaiLiao}
  \textit{
  Suppose that Assumption~\ref{ass.f} holds and that Algorithm~\ref{alg.lmsd} is run with history length $m = 1$ and termination tolerance $\epsilon = 0$.  Then, either $g_{k,1} = 0$ for some $k \in \N{}$ or
  \bequationn
    \|g_{k,1}\| \leq c_1c_2^k\|g_{1,1}\|\ \ \text{for all}\ \ k \in \N{},
  \eequationn
  where, with $K \in \N{}$ from Lemma~\ref{lem.DaiLiao}, the constants are defined as
  \bequationn
    c_1 := 2\(\frac{\lambda_n}{\lambda_1} - 1\)^{K-1}\ \ \text{and}\ \ c_2 := 2^{-1/K} \in (0,1).
  \eequationn
  Overall, the computed gradients vanish $R$-linearly with constants that depend only on $(\lambda_1,\lambda_n)$.
  }
\etheorem

\subsection{$R$-Linear Convergence Rate of LMSD for Arbitrary $m \in \N{}$}

Our goal in this subsection is to build upon the proofs of the results stated in the previous subsection (as given in the cited references) to show that Algorithm~\ref{alg.lmsd} possesses an $R$-linear rate of convergence for any $m \in \N{}$.  More precisely, our goal is to show that the gradients computed by the algorithm vanish $R$-linearly with constants that depend only on the spectrum of the data matrix $A$.  Formally, for simplicity and brevity in our analysis, we make the following standing assumption throughout this section.

\bassumption\label{ass.lmsd}
  \textit{
  Assumption~\ref{ass.f} holds, as do the following:
  \benumerate
    \item[(i)] Algorithm~\ref{alg.lmsd} is run with $\epsilon = 0$ and $g_{k,j} \neq 0$ for all $(k,j) \in \N{} \times [m]$.
    \item[(ii)] For all $k \in \N{}$, the matrix $G_k$ has linearly independent columns.  Further, there exists a scalar~$\rho \geq 1$ such that, for all $k \in \N{}$, the nonsingular matrix $R_k$ satisfies $\|R_k^{-1}\| \leq \rho\|g_{k,1}\|^{-1}$.
  \eenumerate
  }
\eassumption

\noindent
Assumption~\ref{ass.lmsd}$(i)$ is reasonable since, in any situation in which the algorithm terminates finitely, all of our results hold for the iterations prior to that in which the algorithm terminates.  Hence, by proving that the algorithm possesses an $R$-linear rate of convergence for cases when it does not terminate finitely, we claim that it possesses such a rate in all cases.  As for Assumption~\ref{ass.lmsd}$(ii)$, first recall Remark~\ref{rem.G_lin_ind}.  In addition, the bound on the norm of the inverse of $R_k$ is reasonable since, in the case of $m=1$, one finds that $Q_kR_k = G_k = g_{k,1}$ has $Q_k = g_{k,1}/\|g_{k,1}\|$ and $R_k = \|g_{k,1}\|$, meaning that the bound holds with~$\rho=1$.  (This means that, in practice, one might choose $\rho \geq 1$ and iteratively remove columns of $G_k$ for the computation of $T_k$ until one finds $\|R_k^{-1}\| \leq \rho\|g_{k,1}\|^{-1}$, knowing that, in the extreme case, there will remain one column for which this condition is satisfied.  However, for the reasons already given in Remark~\ref{rem.G_lin_ind}, we make Assumption~\ref{ass.lmsd}, meaning that $G_k$ always has $m$ columns.)

We begin by stating two results that reveal important properties of the eigenvalues (corresponding to the elements of $\{T_k\}$) computed by the algorithm, which in turn reveal properties of the stepsizes.  The first result is a direct consequence of the \emph{Cauchy Interlacing Theorem}.  Since this theorem is well-known---see, e.g., \cite{Parl98}---we state the lemma without proof.

\blemma\label{lem.interlacing}
  \textit{
  For all $k \in \N{}$, the eigenvalues of $T_k$ ($= Q_k^TAQ_k$ where $Q_k^TQ_k = I$) satisfy
  \bequationn
    \theta_{k,j} \in [\lambda_{m+1-j}, \lambda_{n+1-j}]\ \ \text{for all}\ \ j \in [m].
  \eequationn
  }
\elemma

The second result provides more detail about how the eigenvalues computed by the algorithm at the end of iteration $k \in \N{}$ relate to the weights in \eqref{eq.combination} corresponding to $k$ for all $j \in [m]$.

\blemma\label{lem.eigenvectors}
  \textit{
  For all $(k,j) \in \N{} \times [m]$, let $q_{k,j} \in \R{m}$ denote the unit eigenvector corresponding to the eigenvalue $\theta_{k,j}$ of $T_k$, i.e., the vector satisfying $T_kq_{k,j} = \theta_{k,j}q_{k,j}$ and $\|q_{k,j}\| = 1$.  Then, defining
  \bequation\label{eq.Dc}
    D_k := \bbmatrix d_{k,1,1} & \cdots & d_{k,m,1} \\ \vdots & \ddots & \vdots \\ d_{k,1,n} & \cdots & d_{k,m,n} \ebmatrix\ \ \text{and}\ \ c_{k,j} := D_kR_k^{-1}q_{k,j},
  \eequation
  it follows that, with the diagonal matrix of eigenvalues (namely, $\Lambda$) defined in Assumption~\ref{ass.f},
  \bequation\label{eq.theta}
    \theta_{k,j} = c_{k,j}^T\Lambda c_{k,j}\ \ \text{and}\ \ c_{k,j}^Tc_{k,j} = 1.
  \eequation
  }
\elemma
\bproof
  For any $k \in \N{}$, it follows from \eqref{eq.Dc} and Lemma~\ref{lem.recursion} (in particular, \eqref{eq.recursion}) that $G_k = QD_k$ where~$Q$ is the orthogonal matrix defined in Assumption~\ref{ass.f}.  Then, since $G_k = Q_kR_k$ (recall \eqref{eq.G}), it follows that $Q_k = QD_kR_k^{-1}$, according to which one finds
  \bequationn
    T_k = Q_k^TAQ_k = R_k^{-T}D_k^TQ^TAQD_kR_k^{-1} = R_k^{-T}D_k^T\Lambda D_kR_k^{-1}.
  \eequationn
  Hence, for each $j \in [m]$, the first equation in \eqref{eq.theta} follows since
  \bequationn
    \theta_{k,j} = q_{k,j}^TT_kq_{k,j} = q_{k,j}^TR_k^{-T}D_k^T\Lambda D_kR_k^{-1}q_{k,j} = c_{k,j}^T\Lambda c_{k,j}.
  \eequationn
  In addition, since $G_k = QD_k$ and the orthogonality of $Q$ imply that $D_k^TD_k = G_k^TG_k$, and since $Q_k = G_kR_k^{-1}$ with $Q_k$ having orthonormal columns (i.e., with $Q_k$ satisfying $Q_k^TQ_k = I$), it follows that
  \bequationn
    c_{k,j}^Tc_{k,j} = q_{k,j}^TR_k^{-T}D_k^TD_kR_k^{-1}q_{k,j} = q_{k,j}^TR_k^{-T}G_k^TG_kR_k^{-1}q_{k,j} = q_{k,j}^TQ_k^TQ_kq_{k,j} = q_{k,j}^Tq_{k,j} = 1,
  \eequationn
  which yields the second equation in \eqref{eq.theta}.
\eproof

The implications of Lemma~\ref{lem.eigenvectors} are seen later in our analysis.  For now, combining Lemma~\ref{lem.interlacing}, Lemma~\ref{lem.recursion} (in particular, \eqref{eq.recursion}), and the fact that \eqref{eq.combination} implies
\bequation\label{eq.g_norm}
  \|g_{k,j}\|^2 = \sum_{i=1}^n d_{k,j,i}^2\ \ \text{for all}\ \ (k,j) \in \N{} \times [m],
\eequation
one is lead to the following result pertaining to recursive properties of the weights in \eqref{eq.combination}.

\blemma\label{lem.loose_bounds}
  \textit{
  For each $(k,j,i) \in \N{} \times [m] \times [n]$, it follows that
  \bequation\label{eq.loose_bound}
    |d_{k,j+1,i}| \leq \delta_{j,i}|d_{k,j,i}|\ \ \text{where}\ \ \delta_{j,i} := \max\left\{\left| 1 - \frac{\lambda_i}{\lambda_{m+1-j}} \right|, \left| 1 - \frac{\lambda_i}{\lambda_{n+1-j}} \right| \right\}.
  \eequation
  Hence, for each $(k,j,i) \in \N{} \times [m] \times [n]$, it follows that
  \bequation\label{eq.loose_bound_m}
    |d_{k+1,j,i}| \leq \Delta_i|d_{k,j,i}|\ \ \text{where}\ \ \Delta_i := \prod_{j=1}^m \delta_{j,i}.
  \eequation
  Furthermore, for each $(k,j,p) \in \N{} \times [m] \times [n]$, it follows that
  \bequation\label{eq.loose_bound_j}
    \sqrt{\sum_{i=1}^p d_{k,j+1,i}^2} \leq \hat\delta_{j,p}\sqrt{\sum_{i=1}^p d_{k,j,i}^2}\ \ \text{where}\ \ \hat\delta_{j,p} := \max_{i\in[p]} \delta_{j,i},
  \eequation
  while, for each $(k,j) \in \N{} \times [m]$, it follows that
  \bequation\label{eq.loose_bound_g}
    \|g_{k+1,j}\| \leq \Delta \|g_{k,j}\|\ \ \text{where}\ \ \Delta := \max_{i\in[n]} \Delta_i.
  \eequation
  }
\elemma
\bproof
  Recall that, for any given $(k,j,i) \in \N{} \times [m] \times [n]$, Lemma~\ref{lem.recursion} (in particular,  \eqref{eq.recursion}) states
  \bequationn
    d_{k,j+1,i} = (1 - \alpha_{k,j}\lambda_i)d_{k,j,i}.
  \eequationn
  The relationship \eqref{eq.loose_bound} then follows due to Lemma~\ref{lem.interlacing}, which, in particular, shows that
  \bequationn
    \alpha_{k,j} \in \left[\frac{1}{\lambda_{n+1-j}},\frac{1}{\lambda_{m+1-j}}\right] \subseteq \left[\frac{1}{\lambda_n},\frac{1}{\lambda_1}\right]\ \ \text{for all}\ \ (k,j) \in \N{} \times [m].
  \eequationn
  The consequence \eqref{eq.loose_bound_m} then follows by combining \eqref{eq.loose_bound} for all $j \in [m]$ and recalling that Step~\ref{step.end_of_cycle} yields $g_{k+1,1} \gets g_{k,m+1}$ for all $k \in \N{}$.  Now, from \eqref{eq.loose_bound}, one finds that
  \bequationn
    \sum_{i=1}^p d_{k,j+1,i}^2 \leq \sum_{i=1}^p \delta_{j,i}^2 d_{k,j,i}^2 \leq \hat\delta_{j,p}^2 \sum_{i=1}^p d_{k,j,i}^2\ \ \text{for all}\ \ (k,j,p) \in \N{} \times [m] \times [n],
  \eequationn
  yielding the desired conclusion \eqref{eq.loose_bound_j}.  Finally, combining \eqref{eq.loose_bound_m} and \eqref{eq.g_norm}, one obtains that
  \bequationn
    \|g_{k+1,j}\|^2 = \sum_{i=1}^n d_{k+1,j,i}^2 \leq \sum_{i=1}^n \Delta_i^2 d_{k,j,i}^2 \leq \Delta^2 \sum_{i=1}^n d_{k,j,i}^2 = \Delta^2 \|g_{k,j}\|^2\ \ \text{for all}\ \ (k,j) \in \N{} \times [m],
  \eequationn
  yielding the desired conclusion \eqref{eq.loose_bound_g}.
\eproof

A consequence of the previous lemma is that if $\Delta_i \in [0,1)$ for all $i \in [n]$, then $\Delta \in [0,1)$, from which~\eqref{eq.loose_bound_g} implies that, for each $j \in [m]$, the gradient norm sequence $\{\|g_{k,j}\|\}_{k\in\N{}}$ vanishes $Q$-linearly.  For example, such a situation occurs when $\lambda_n < 2\lambda_1$.  However, as noted in \cite{DaiLiao02}, this is a highly special case that should not be assumed to hold widely in practice.  A more interesting and widely relevant consequence of the lemma is that for any $i \in [n]$ such that $\Delta_i \in [0,1)$, the sequences $\{|d_{k,j,i}|\}_{k\in\N{}}$ for each $j \in [m]$ vanish $Q$-linearly.  For example, this is \emph{always} true for $i=1$, where
\bequationn
  \delta_{j,1} = \max\left\{1 - \frac{\lambda_1}{\lambda_{m+1-j}},1 - \frac{\lambda_1}{\lambda_{n+1-j}}\right\} \in [0,1)\ \ \text{for all}\ \ j \in [m],
\eequationn
from which it follows that
\bequation\label{eq.Delta_1}
  \Delta_1 = \prod_{j=1}^m \delta_{j,1} \in [0,1).
\eequation
The following is a crucial consequence that one can draw from this observation.

\blemma\label{lem.i=1}
  \textit{
  If $\Delta_1 = 0$, then $d_{1+\khat,\jhat,1} = 0$ for all $(\khat,\jhat) \in \N{} \times [m]$.  Otherwise, if $\Delta_1 > 0$, then:
  \benumerate
    \item[(i)] for any $(k,j) \in \N{} \times [m]$ such that $d_{k,j,1} = 0$, it follows that $d_{k+\khat,\jhat,1} = 0$ for all $(\khat,\jhat) \in \N{} \times [m]$;
    \item[(ii)] for any $(k,j) \in \N{} \times [m]$ such that $|d_{k,j,1}| > 0$ and any $\epsilon_1 \in (0,1)$, it follows that
    \bequationn
      \frac{|d_{k+\khat,\jhat,1}|}{|d_{k,j,1}|} \leq \epsilon_1\ \ \text{for all}\ \ \khat \geq 1 + \left\lceil \frac{\log\epsilon_1}{\log\Delta_1} \right\rceil\ \ \text{and}\ \ \jhat \in [m].
    \eequationn
  \eenumerate
  }
\elemma
\bproof
  If $\Delta_1 = 0$, then the desired conclusion follows from Lemma~\ref{lem.loose_bounds}; in particular, it follows from the inequality~\eqref{eq.loose_bound_m} for $i = 1$.  Similarly, for any $(k,j) \in \N{} \times [m]$ such that $d_{k,j,1} = 0$, the conclusion in part~$(i)$ follows from the same conclusion in Lemma~\ref{lem.loose_bounds}, namely, \eqref{eq.loose_bound_m} for $i=1$.  Hence, let us continue to prove part $(ii)$ under the assumption that $\Delta_1 \in (0,1)$ (recall \eqref{eq.Delta_1}).

  Suppose that the given condition holds with $j=1$, i.e., consider $k \in \N{}$ such that $|d_{k,1,1}| > 0$.  Then, it follows by Lemma~\ref{lem.loose_bounds} (in particular, \eqref{eq.loose_bound_m} for $j=1$ and $i=1$) that
  \bequation\label{eq.Delta}
    \frac{|d_{k+\khat,1,1}|}{|d_{k,1,1}|} \leq \Delta_1^{\khat}\ \ \text{for any}\ \ \khat \in \N{}.
  \eequation
  Since $\Delta_1 \in (0,1)$, taking the logarithm of the term on the right-hand side with $\khat = \lceil \log\epsilon_1/\log\Delta_1 \rceil$ yields
  \bequation\label{eq.log}
    \left\lceil \frac{\log\epsilon_1}{\log\Delta_1} \right\rceil \log\Delta_1 \leq \(\frac{\log\epsilon_1}{\log\Delta_1}\)\log\Delta_1 = \log\(\epsilon_1\).
  \eequation
  Since $\log(\cdot)$ is nondecreasing, the inequalities yielded by \eqref{eq.log} combined with \eqref{eq.Delta} along with \eqref{eq.loose_bound_m} from Lemma~\ref{lem.loose_bounds} yield the desired result for $j=1$.  On the other hand, if the conditions of part $(ii)$ hold for some other $j \in [m]$, then the desired conclusion follows from a similar reasoning, though an extra cycle may need to be completed before the desired conclusion holds for all points in the cycle, i.e., for all $\jhat \in [m]$; hence the addition of 1 to $\lceil \log\epsilon_1/\log\Delta_1 \rceil$ in the general conclusion.
\eproof

One may conclude from Lemma~\ref{lem.i=1} and \eqref{eq.combination} that, for any $(k,j) \in \N{} \times [m]$ and $\epsilon_1 \in (0,1)$, one has
\bequationn
  \frac{|d_{k+\khat,\jhat,1}|}{\|g_{k,j}\|} \leq \epsilon_1\ \ \text{for all}\ \ \khat \geq K_1\ \ \text{and}\ \ \jhat \in [m]
\eequationn
for some $K_1 \in \N{}$ that depends on the desired contraction factor $\epsilon_1 \in (0,1)$ and the problem-dependent constant~$\Delta_1 \in (0,1)$, but does \emph{not} depend on the iteration number pair $(k,j)$.  Our goal now is to show that if a similar, but looser conclusion holds for a squared sum of the weights in~\eqref{eq.combination} up through $p \in [n-1]$, then the squared weight corresponding to index $p + 1$ eventually becomes sufficiently small in a number of iterations that is independent of the iteration number $k$.  (For this lemma, we fix $j=\jhat=1$ so as to consider only the first gradient in each cycle.  This choice is somewhat arbitrary since our concluding theorem will confirm that a similar result would hold for any $j \in [m]$ and $\jhat = j$.)  For the lemma, we define the following constants that dependent only on $p$, the spectrum of $A$ (which, in particular, yields the bounds and definitions in \eqref{lem.loose_bounds}), and the scalar constant $\rho \geq 1$ from Assumption~\ref{ass.lmsd}:
\bsubequations\label{eq.get_small_defs}
  \begin{align}
    \hat\delta_p &:= \(1 + \hat\delta_{1,p}^2 + \hat\delta_{1,p}^2\hat\delta_{2,p}^2 + \cdots + \prod_{j=1}^{m-1} \hat\delta_{j,p}^2 \) \in [1,\infty), \\
    \hat\Delta_{p+1} &:= \max\left\{\frac13,1 - \frac{\lambda_{p+1}}{\lambda_n}\right\}^m \in (0,1), \label{eq.Deltahatp1} \\
    \text{and}\ \ \hat{K}_p &:= \left\lceil \frac{\log\(2\hat\delta_p\rho\epsilon_p\Delta_{p+1}^{-(K_p+1)}\)}{\log \hat\Delta_{p+1}} \right\rceil. \label{def.Khatp}
  \end{align}
\esubequations

\blemma\label{lem.get_small}
  \textit{
  For any $(k,p) \in \N{} \times [n-1]$, if there exists $(\epsilon_p,K_p) \in (0,\tfrac{1}{2\hat\delta_p\rho}) \times \N{}$ independent of $k$ with
  \bequation\label{eq.p}
    \sum_{i=1}^p d_{k+\khat,1,i}^2 \leq \epsilon_p^2 \|g_{k,1}\|^2\ \ \text{for all}\ \ \khat \geq K_p,
  \eequation
  then one of the following holds:
  \benumerate
    \item[(i)] $\Delta_{p+1} \in [0,1)$ and there exists $K_{p+1} \geq K_p$ dependent only on $\epsilon_p$, $\rho$, and the spectrum of $A$ with
    \bequation\label{eq.p+1_easy}
      d_{k+K_{p+1},1,p+1}^2 \leq 4\hat\delta_p^2\rho^2\epsilon_p^2\|g_{k,1}\|^2;
    \eequation
    \item[(ii)] $\Delta_{p+1} \in [1,\infty)$ and, with $K_{p+1} := K_p + \hat{K}_p + 1$, there exists $\khat_0 \in \{K_p,\dots,K_{p+1}\}$ with
    \bequation\label{eq.p+1}
      d_{k+\khat_0,1,p+1}^2 \leq 4\hat\delta_p^2\rho^2\epsilon_p^2\|g_{k,1}\|^2.
    \eequation
  \eenumerate
  }
\elemma
\bproof
  By Lemma~\ref{lem.loose_bounds} (in particular, \eqref{eq.loose_bound_m} with $j=1$ and $i = p+1$) and \eqref{eq.g_norm}, it follows that
  \bequation\label{eq.jump}
    d_{k+\khat,1,p+1}^2 \leq \(\Delta_{p+1}^{\khat} d_{k,1,p+1}\)^2 = \Delta_{p+1}^{2\khat} d_{k,1,p+1}^2 \leq \Delta_{p+1}^{2\khat} \|g_{k,1}\|^2\ \ \text{for all}\ \ \khat \in \N{}.
  \eequation
  If $\Delta_{p+1} \in [0,1)$, then \eqref{eq.jump} immediately implies the existence of $K_{p+1}$ dependent only on $\epsilon_p$, $\rho$, and the spectrum of $A$ such that \eqref{eq.p+1_easy} holds.  Hence, let us continue under the assumption that $\Delta_{p+1} \geq 1$, where one should observe that $\rho \geq 1$, $\hat\delta_p \geq 1$, $\epsilon_p \in (0,\tfrac{1}{2\hat\delta_p\rho})$, $K_p \in \N{}$, and $\Delta_{p+1} \geq 1$ imply $2\hat\delta_p\rho\epsilon_p\Delta_{p+1}^{-K_p} \in (0,1)$, meaning that $\hat{K}_p \in \N{}$.  To prove the desired result, it suffices to show that if
  \bequation\label{eq.stay_big}
    d_{k+\khat,1,p+1}^2 > 4\hat\delta_p^2\rho^2\epsilon_p^2\|g_{k,1}\|^2\ \ \text{for all}\ \ \khat \in \{K_p,\dots,K_{p+1}-1\},
  \eequation
  then \eqref{eq.p+1} holds at the beginning of the next cycle (i.e., when $\khat_0 = K_{p+1}$).  From Lemma~\ref{lem.eigenvectors}, Lemma~\ref{lem.loose_bounds} (in particular, \eqref{eq.loose_bound_j}), \eqref{eq.p}, and \eqref{eq.stay_big}, it follows that with $\{c_{k+\khat,j,i}\}_{i=1}^n$ representing the elements of the vector $c_{k+\khat,j}$ and the matrix $D_{k+\khat,p}$ representing the first $p$ rows of $D_{k+\khat}$, one finds
  \bequalin
    \sum_{i=1}^p c_{k+\khat,j,i}^2
      &\leq \|D_{k+\khat,p}\|_2^2 \|R_{k+\khat}^{-1}\|^2\|q_{k+\khat,j}\|^2 \\
      &\leq \(1 + \hat\delta_{1,p}^2 + \hat\delta_{1,p}^2\hat\delta_{2,p}^2 + \cdots + \prod_{j=1}^{m-1} \hat\delta_{j,p}^2 \) \(\sum_{i=1}^p d_{k+\khat,1,i}^2\) \rho^2 \|g_{k+\khat,1}\|^{-2} \\
      &\leq \hat\delta_p^2 (\epsilon_p^2 \|g_{k,1}\|^2) \rho^2 (4\hat\delta_p^2\rho^2\epsilon_p^2)^{-1}\|g_{k,1}\|^{-2} \leq \tfrac14\ \ \text{for all}\ \ \khat \in \{K_p,\dots,K_{p+1}-1\}\ \ \text{and}\ \ j \in [m].
  \eequalin
  Along with Lemma~\ref{lem.eigenvectors}, this implies that
  \bequation\label{eq.3/4}
    \theta_{k+\khat,j} = \sum_{i=1}^n \lambda_i c_{k+\khat,j,1}^2 \geq \tfrac34 \lambda_{p+1}\ \ \text{for all}\ \ \khat \in \{K_p,\dots,K_{p+1}-1\}\ \ \text{and}\ \ j \in [m].
  \eequation
  Together with Lemma~\ref{lem.recursion} (see \eqref{eq.recursion}) and $\alpha_{k+\khat+1,j} = \theta_{k+\khat,j}^{-1}$ for all $j \in [m]$, the bound \eqref{eq.3/4} implies
  \begin{align}
    d_{k+\khat+2,1,p+1}^2
      &= \(\prod_{j=1}^m \(1 - \alpha_{k+\khat+1,j}\lambda_{p+1}\)^2\) d_{k+\khat+1,1,p+1}^2 \nonumber \\
      &\leq \hat\Delta_{p+1}^2 d_{k+\khat+1,1,p+1}^2 \ \ \text{for all}\ \ \khat \in \{K_p,\dots,K_{p+1}-1\}. \label{eq.scream}
  \end{align}
  Applying this bound recursively, it follows with $K_{p+1} = K_p + \hat{K}_p + 1$ and \eqref{eq.jump} for $\khat = K_{p+1}$ that
  \bequationn
    d_{k+K_{p+1},1,p+1}^2 \leq \hat\Delta_{p+1}^{2\hat{K}_p} d_{k+K_p+1,1,p+1}^2 \leq \hat\Delta_{p+1}^{2\hat{K}_p} \Delta_{p+1}^{2(K_p+1)} \|g_{k,1}\|^2 \leq 4\hat\delta_p^2r^2\epsilon_p^2\|g_{k,1}\|^2,
  \eequationn
  where the last inequality follows by the definition of $\Khat_p$ in \eqref{def.Khatp}.
\eproof

We have shown that small squared weights in~\eqref{eq.combination} associated with indices up through $p \in [n-1]$ imply that the squared weight associated with index $p+1$ eventually becomes small.  The next lemma shows that these latter squared weights also remain sufficiently small indefinitely.

\blemma\label{lem.stay_small}
  \textit{
  For any $(k,p) \in \N{} \times [n-1]$, if there exists $(\epsilon_p,K_p) \in (0,\tfrac{1}{2\hat\delta_p\rho}) \times \N{}$ independent of $k$ such that \eqref{eq.p} holds, then, with $\epsilon_{p+1}^2 := (1 + 4\max\{1,\Delta_{p+1}^4\}\hat\delta_p^2\rho^2)\epsilon_p^2$ and $K_{p+1} \in \N{}$ from Lemma~\ref{lem.get_small},
  \bequation\label{eq.p+1_all}
    \sum_{i=1}^{p+1} d_{k+\khat,1,i}^2 \leq \epsilon_{p+1}^2\|g_{k,1}\|^2\ \ \text{for all}\ \ \khat \geq K_{p+1}.
  \eequation
  }
\elemma
\bproof
  For the same reasons as in the proof of Lemma~\ref{lem.get_small}, the result follows if $\Delta_{p+1} \in [0,1)$.  Hence, we may continue under the assumption that $\Delta_{p+1} \geq 1$ and define $\hat\Delta_{p+1} \in (0,1)$ and $\Khat_p \in \N{}$ as in~\eqref{eq.get_small_defs}.  By Lemma~\ref{lem.get_small}, there exists $\khat_0 \in \{K_p,\dots,K_{p+1}\}$ such that
  \bequation\label{eq.p+1_k0}
    d_{k+\khat,1,p+1}^2 \leq 4\hat\delta_p^2\rho^2\epsilon_p^2\|g_{k,1}\|^2\ \ \text{when}\ \ \khat = \khat_0.
  \eequation
  If the inequality in \eqref{eq.p+1_k0} holds for all $\khat \geq \khat_0$, then \eqref{eq.p+1_all} holds with $\epsilon_{p+1}^2 = (1 + 4\hat\delta_p^2\rho^2)\epsilon_p^2$.  Otherwise, let $\khat_1 \in \N{}$ denote the smallest natural number such that
  \bequation\label{eq.shoulder}
    d_{k+\khat,1,p+1}^2 \leq 4\hat\delta_p^2\rho^2\epsilon_p^2\|g_{k,1}\|^2\ \ \text{for all}\ \ \khat_0 \leq \khat \leq \khat_1,
  \eequation
  but
  \bequation\label{eq.large_again}
    d_{k+\khat_1+1,1,p+1}^2 > 4\hat\delta_p^2\rho^2\epsilon_p^2\|g_{k,1}\|^2.
  \eequation
  As in the arguments that lead to \eqref{eq.scream} in the proof of Lemma~\ref{lem.get_small}, combining \eqref{eq.p} and \eqref{eq.large_again} implies
  \bequationn
    d_{k+\khat_1+3,1,p+1}^2 \leq \hat\Delta_{p+1}^2 d_{k+\khat_1+2,1,p+1}^2.
  \eequationn
  Generally, this same argument can be used to show that
  \bequationn
    \khat \geq K_p\ \ \text{and}\ \ d_{k+\khat+1,1,p+1}^2 > 4\hat\delta_p^2\rho^2\epsilon_p^2\|g_{k,1}\|^2\ \ \text{imply}\ \ d_{k+\khat+3,1,p+1}^2 \leq \hat\Delta_{p+1}^2 d_{k+\khat+2,1,p+1}^2.
  \eequationn
  Since $\hat\Delta_{p+1} \in (0,1)$, this fact and \eqref{eq.large_again} imply the existence of $\khat_2 \in \N{}$ such that
  \bequation\label{eq.airplane}
    d_{k+\khat+1,1,p+1}^2 > 4\hat\delta_p^2\rho^2\epsilon_p^2\|g_{k,1}\|^2\ \ \text{for all}\ \ \khat_1 \leq \khat \leq \khat_2 - 2,
  \eequation
  but
  \bequationn
    d_{k+\khat_2,1,p+1}^2 \leq 4\hat\delta_p^2\rho^2\epsilon_p^2\|g_{k,1}\|^2,
  \eequationn
  while, from above,
  \bequation\label{eq.ground}
    d_{k+\khat+3,1,p+1}^2 \leq \hat\Delta_{p+1}^2 d_{k+\khat+2,1,p+1}^2\ \ \text{for all}\ \ \khat_1 \leq \khat \leq \khat_2 - 2.
  \eequation
  Moreover, by Lemma~\ref{lem.loose_bounds} (in particular, \eqref{eq.loose_bound_m}) and \eqref{eq.shoulder}, it follows that
  \bsubequations
    \begin{align}
      d_{k+\khat_1+1,1,p+1}^2 &\leq \Delta_{p+1}^2 d_{k+\khat_1,1,p+1}^2 \leq 4\Delta_{p+1}^2\hat\delta_p^2\rho^2\epsilon_p^2\|g_{k,1}\|^2 \\
      \text{and}\ \ d_{k+\khat_1+2,1,p+1}^2 &\leq 4\Delta_{p+1}^4\hat\delta_p^2\rho^2\epsilon_p^2\|g_{k,1}\|^2. \label{eq.plastic}
    \end{align}
  \esubequations
  Combining \eqref{eq.ground} and \eqref{eq.plastic}, it follows that
  \bequationn
    d_{k+\khat+3,1,p+1}^2 \leq 4\hat\Delta_{p+1}^2\Delta_{p+1}^4\hat\delta_p^2\rho^2\epsilon_p^2\|g_{k,1}\|^2\ \ \text{for all}\ \ \khat_1 \leq \khat \leq \khat_2 - 2.
  \eequationn
  Overall, since \eqref{eq.Deltahatp1} ensures $\hat\Delta_{p+1} \in (0,1)$, we have shown that
  \bequation\label{eq.suffice}
    d_{k+\khat,1,p+1}^2 \leq 4\Delta_{p+1}^4\hat\delta_p^2\rho^2\epsilon_p^2 \|g_{k,1}\|^2\ \ \text{for all}\ \ \khat \in \{\khat_0,\dots,\khat_2\}.
  \eequation
  Repeating this argument for later iterations, we arrive at the desired conclusion.
\eproof

The following lemma is a generalization of Lemma~\ref{lem.DaiLiao} for any $m \in \N{}$.  Our proof is similar to that of Lemma~2.4 in \cite{DaiLiao02}.  We provide it in full for completeness.

\blemma\label{lem.g_contraction}
  \textit{
  There exists $K \in \N{}$ dependent only on the spectrum of $A$ such that
  \bequationn
    \|g_{k+K,1}\| \leq \thalf \|g_{k,1}\|\ \ \text{for all}\ \ k \in \N{}.
  \eequationn
  }
\elemma
\bproof
  By Lemma~\ref{lem.stay_small}, if for some $(\epsilon_p,K_p) \in (0,\tfrac{1}{2\hat\delta_p\rho}) \times \N{}$ independent of $k$ one finds
  \bequation\label{eq.p_again}
    \sum_{i=1}^p d_{k+\khat,1,i}^2 \leq \epsilon_p^2 \|g_{k,1}\|^2\ \ \text{for all}\ \ \khat \geq K_p,
  \eequation
  then for $\epsilon_{p+1}^2 := (1 + 4\max\{1,\Delta_{p+1}^4\}\hat\delta_p^2\rho^2)\epsilon_p^2$ and some $K_{p+1} \geq K_p$ independent of $k$ one finds
  \bequation\label{eq.p+1_all_again}
    \sum_{i=1}^{p+1} d_{k+\khat,1,i}^2 \leq \epsilon_{p+1}^2\|g_{k,1}\|^2\ \ \text{for all}\ \ \khat \geq K_{p+1}.
  \eequation
  Since Lemma~\ref{lem.i=1} implies that for any $\epsilon_1 \in (0,1)$ one can find $K_1$ independent of $k$ such that \eqref{eq.p_again} holds with $p=1$, it follows that, independent of $k$, there exists a sufficiently small $\epsilon_1 \in (0,1)$ such that
  \bequationn
    \epsilon_1^2 \leq \cdots \leq \epsilon_n^2 \leq \tfrac14.
  \eequationn
  Hence, for any $k \in \N{}$, it follows that there exists $K = K_n$ such that
  \bequationn
    \|g_{k+\khat,1}\|^2 = \sum_{i=1}^n d_{k+\khat,1,i}^2 \leq \tfrac14 \|g_{k,1}\|^2\ \ \text{for all}\ \ \khat \geq K,
  \eequationn
  as desired.
\eproof

We are now prepared to state our final result, the proof of which follows in the same manner as Theorem~\ref{th.DaiLiao} follows from Lemma~\ref{lem.DaiLiao} in \cite{DaiLiao02}.  We prove it in full for completeness.

\btheorem
  \textit{
  The sequence $\{\|g_{k,1}\|\}$ vanishes $R$-linearly.
  }
\etheorem
\bproof
  If $\Delta \in [0,1)$, then it has already been argued (see the discussion following Lemma~\ref{lem.loose_bounds}) that~$\{\|g_{k,1}\|\}$ vanishes $Q$-linearly.  Hence, let us continue assuming that $\Delta \geq 1$.  By Lemma~\ref{lem.g_contraction}, there exists~$K \in \N{}$ dependent only on the spectrum of $A$ such that
  \bequationn
    \|g_{1+Kl,1}\| \leq \thalf\|g_{1+K(l-1),1}\|\ \ \text{for all}\ \ l \in \N{}.
  \eequationn
  Applying this result recursively, it follows that
  \bequation\label{eq.tight_bound_g}
    \|g_{1+Kl,1}\| \leq (\thalf)^l \|g_{1,1}\|\ \ \text{for all}\ \ l \in \N{}.
  \eequation
  Now, for any $k \geq 1$, let us write $k = Kl+\khat$ for some $l \in \{0\}\cup\N{}$ and $\khat \in \{0\}\cup[K-1]$.  It follows that
  \bequationn
    l = k/K - \khat/K \geq k/K - 1.
  \eequationn
  By this fact, \eqref{eq.loose_bound_g}, and \eqref{eq.tight_bound_g}, it follows that for any $k = Kl+\khat \in \N{}$ one has
  \bequationn
    \|g_{k,1}\| \leq \Delta^{\khat-1}\|g_{1+Kl,1}\| \leq \Delta^{K-1}(\thalf)^{k/K-1}\|g_{1,1}\| \leq c_1c_2^k\|g_{1,1}\|,
  \eequationn
  where
  \bequationn
    c_1 := 2\Delta^{K-1} \ \ \text{and}\ \ c_2 := 2^{-1/K} \in (0,1),
  \eequationn
  which implies the desired conclusion.
\eproof

\section{Numerical Demonstrations}\label{sec.numerical}

The analysis in the previous section provides additional insights into the behavior of Algorithm~\ref{alg.lmsd} beyond its $R$-linear rate of convergence.  In this section, we provide the results of numerical experiments to demonstrate the behavior of the algorithm in a few types of cases.  The algorithm was implemented and the experiments were performed in Matlab.  It is not our goal to show the performance of Algorithm~\ref{alg.lmsd} for various values of $m$, say to argue whether the performance improves or not as $m$ is increased.  This is an important question for which some interesting discussion is provided by \cite{Flet12}.  However, to determine what is a good choice of $m$ for various types of cases would require a larger set of experiments that are outside of the scope of this paper.  For our purposes, our only goal is to provide some simple illustrations of the behavior as shown by our theoretical analysis.

Our analysis reveals that the convergence behavior of the algorithm depends on the spectrum of the matrix $A$.  Therefore, we have constructed five test examples, all with $n=100$, but with different eigenvalue distributions.  For the first problem, the eigenvalues of $A$ are evenly distributed in $[1,1.9]$.  Since this ensures that $\lambda_n < 2\lambda_1$, our analysis reveals that the algorithm converges $Q$-linearly for this problem; recall the discussion after Lemma~\ref{lem.loose_bounds}.  All other problems were constructed so that $\lambda_1 = 1$ and $\lambda_n = 100$, for which one clearly finds $\lambda_n > 2\lambda_1$.  For the second problem, all eigenvalues are evenly distributed in $[\lambda_1,\lambda_n]$; for the third problem, the eigenvalues are clustered in five distinct blocks; for the fourth problem, all eigenvalues except one are clustered around $\lambda_1$; and for the fifth problem, all eigenvalues except one are clustered around $\lambda_n$.  Table~\ref{tab.test_problems} shows the spectrum of $A$ for each problem.

The table also shows the numbers of outer and (total) inner iterations required by Algorithm~\ref{alg.lmsd} (indicated by column headers ``$k$'' and ``$j$'', respectively) when it was run with $\epsilon = 10^{-8}$ and either $m=1$ or $m=5$.  In all cases, the initial $m$ stepsizes were generated randomly from a uniform distribution over the interval $[\lambda_{100}^{-1},\lambda_1^{-1}]$.  One finds that the algorithm terminates in relatively few outer and inner iterations relative to $n$, especially when many of the eigenvalues are clustered.  This dependence on clustering of the eigenvalues should not be surprising since, recalling Lemma~\ref{lem.interlacing}, clustered eigenvalues makes it likely that an eigenvalue of $T_k$ will be near an eigenvalue of $A$, which in turn implies by Lemma~\ref{lem.recursion} that the weights in the representation \eqref{eq.combination} will vanish quickly.  On the other hand, for the problems for which the eigenvalues are more evenly spread in $[1,100]$, the algorithm requires relatively more outer iterations, though still not an excessively large number relative to $n$.  For these problems, the performance was better for $m=5$ versus $m=1$, both in terms of outer and (total) inner iterations.

\btable[ht]\renewcommand{\tabcolsep}{10pt}
  \centering
  \caption{Spectra of $A$ for five test problems along with outer and (total) inner iteration counts required by Algorithm~\ref{alg.lmsd}.  For each spectrum, a set of eigenvalues in an interval indicates that the eigenvalues are evenly distributed within that interval.}
  \label{tab.test_problems}
  \btabular{|c|rcl|c|c|c|c|}
    \hline
            &                                        &       &                & \multicolumn{2}{c|}{$m=1$} & \multicolumn{2}{c|}{$m=5$} \\
    Problem & \multicolumn{3}{c|}{Spectrum} & \multicolumn{1}{c}{$k$} & $j$ & \multicolumn{1}{c}{$k$} & $j$ \\
    \hline
    \hline
    1       & $\{\lambda_{ 1},\dots,\lambda_{100}\}$ & $\subset$ & $[ 1,  1.9]$ & 13 & 13 & 3 & 14 \\
    \hline
    2       & $\{\lambda_{ 1},\dots,\lambda_{100}\}$ & $\subset$ & $[ 1,100  ]$ & 124 & 124 & 23 & 114 \\
    \hline
    3       & $\{\lambda_{ 1},\dots,\lambda_{ 20}\}$ & $\subset$ & $[ 1,  2  ]$ & 112 & 112 & 16 & 79 \\
            & $\{\lambda_{21},\dots,\lambda_{ 40}\}$ & $\subset$ & $[25, 26  ]$ &  &  &  &  \\
            & $\{\lambda_{41},\dots,\lambda_{ 60}\}$ & $\subset$ & $[50, 51  ]$ &  &  &  &  \\
            & $\{\lambda_{61},\dots,\lambda_{ 80}\}$ & $\subset$ & $[75, 76  ]$ &  &  &  &  \\
            & $\{\lambda_{81},\dots,\lambda_{100}\}$ & $\subset$ & $[99,100  ]$ &  &  &  &  \\
    \hline
    4       & $\{\lambda_{ 1},\dots,\lambda_{ 99}\}$ & $\subset$ & $[ 1,  2  ]$ & 26 & 26 & 4 & 20 \\
            & $\lambda_{100}$ & $=$ & $100$                                   &  &  &  &  \\
    \hline
    5       & $\lambda_1$ & $=$ & $1$                                         & 16 & 16 & 5 & 25 \\
            & $\{\lambda_{ 2},\dots,\lambda_{100}\}$ & $\subset$ & $[99,100  ]$ &  &  &  &  \\
    \hline
  \etabular
\etable

As seen in our analysis (inspired by \cite{Rayd93}, \cite{DaiLiao02}, and \cite{Flet12}), a more refined look into the behavior of the algorithm is obtained by observing the step-by-step magnitudes of the weights in~\eqref{eq.combination} for the generated gradients.  Hence, for each of the test problems, we plot in Figures~\ref{fig.p1}, \ref{fig.p2}, \ref{fig.p3}, \ref{fig.p4}, and \ref{fig.p5} these magnitudes (on a log scale) for a few representative values of $i \in [n]$.  Each figure consists of four sets of plots: the first and third show the magnitudes corresponding to $\{g_{k,1}\}$ (i.e., for the first point in each cycle) when $m=1$ and $m=5$, respectively, while the second and fourth show the magnitudes at all outer and inner iterations, again when $m=1$ and $m=5$, respectively.  In a few of the images, the plot ends before the right-hand edge of the image.  This is due to the log of the absolute value of the weight being evaluated as $-\infty$ in Matlab.

The tables show that the magnitudes of the weights corresponding to $i=1$ always decrease monotonically, as proved in Lemma~\ref{lem.i=1}.  The magnitudes corresponding to $i=2$ also often decrease monotonically, but, as seen in the results for Problem~5, this is not always the case.  In any case, the magnitudes corresponding to $i=50$ and $i=100$ often do not decrease monotonically, though, as proved in our analysis, one observes that the magnitudes demonstrate a downward trend over a finite number of cycles.

Even further insight into the plots of these magnitudes can be gained by observing the values of the constants $\{\Delta_i\}_{i\in[n]}$ for each problem and history length.  Recalling \eqref{eq.loose_bound_m}, these constants bound the increase that a particular weight in \eqref{eq.combination} might experience from one point in a cycle to the same point in the subsequent cycle.  For illustration, we plot in Figures~\ref{fig.p1delta}, \ref{fig.p2delta}, \ref{fig.p3delta}, \ref{fig.p4delta}, and \ref{fig.p5delta} these constants.  Values less than 1 are indicated by a purple bar while values greater than or equal to 1 are indicated by a blue bar.  Note that, in Figure~\ref{fig.p4delta}, all values are small for both history lengths except $\Delta_{100}$.  In Figure~\ref{fig.p5delta}, $\Delta_1$ is less than one in both figures, but the remaining constants are large for $m=1$ while being small for $m=5$.

\bfigure[H]
  \centering
  \includegraphics[width=0.47\textwidth,clip=true,trim=15 5 55 15]{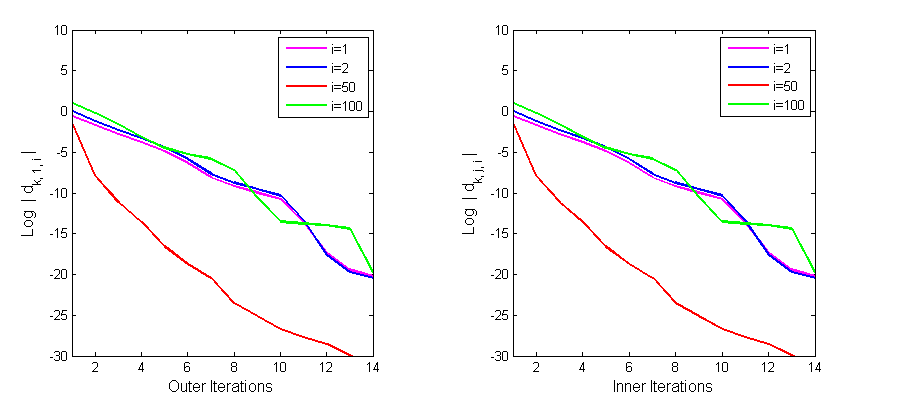}\qquad \includegraphics[width=0.47\textwidth,clip=true,trim=15 5 55 15]{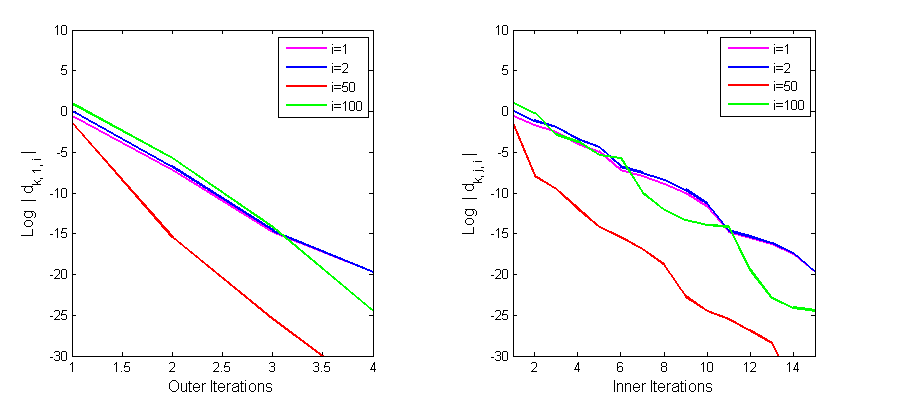}
  \caption{Weights in \eqref{eq.combination} for problem $1$ with history length $m=1$ (left two plots) and $m=5$ (right two plots).}
  \label{fig.p1}
\efigure

\bfigure[H]
  \centering
  \includegraphics[width=0.33\textwidth,clip=true,trim=10 15 25 0]{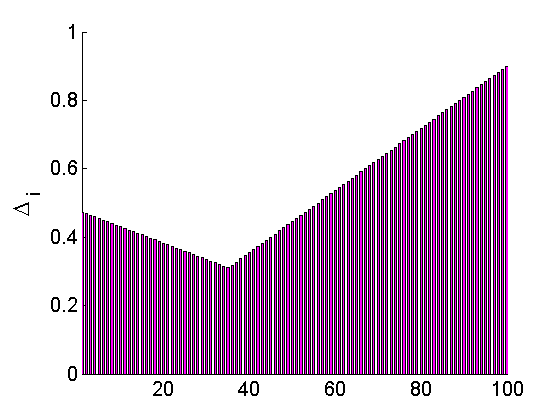}\qquad \includegraphics[width=0.33\textwidth,clip=true,trim=10 15 25 0]{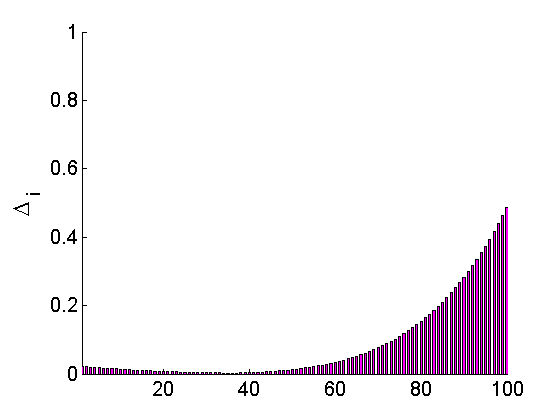}
  \caption{Constants in \eqref{eq.loose_bound_m} for problem $1$ with history length $m=1$ (left plot) and $m=5$ (right plot).}
  \label{fig.p1delta}
\efigure

\bfigure[H]
  \centering
  \includegraphics[width=0.47\textwidth,clip=true,trim=15 5 55 15]{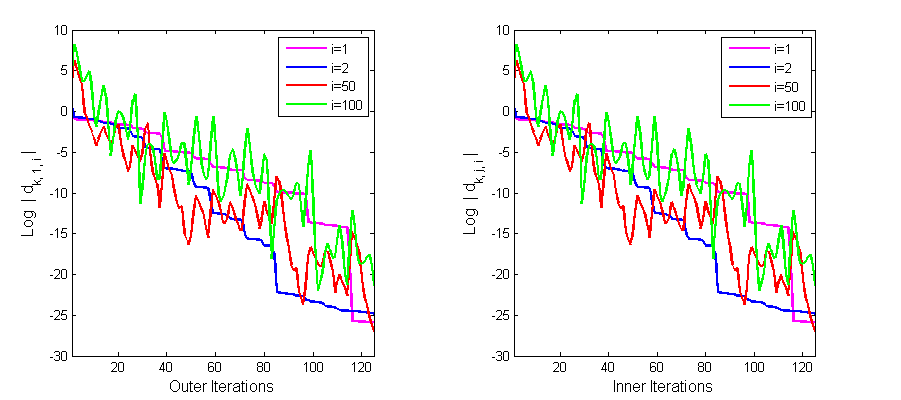}\qquad \includegraphics[width=0.47\textwidth,clip=true,trim=15 5 55 15]{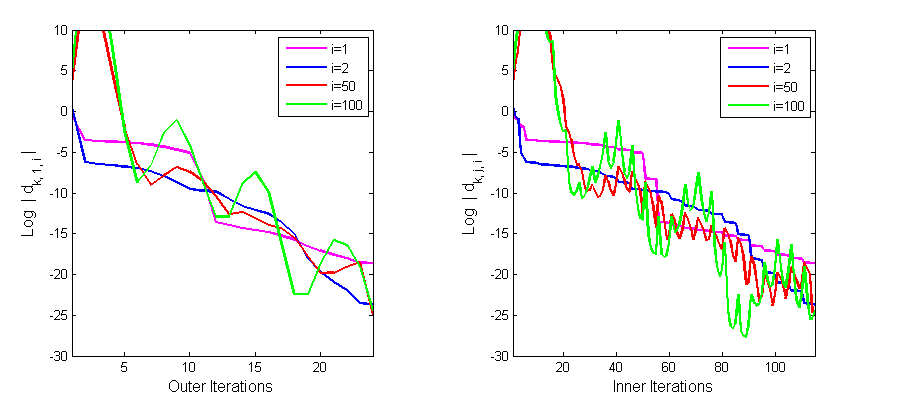}
  \caption{Weights in \eqref{eq.combination} for problem $2$ with history length $m=1$ (left two plots) and $m=5$ (right two plots).}
  \label{fig.p2}
\efigure

\bfigure[H]
  \centering
  \includegraphics[width=0.33\textwidth,clip=true,trim=10 15 25 0]{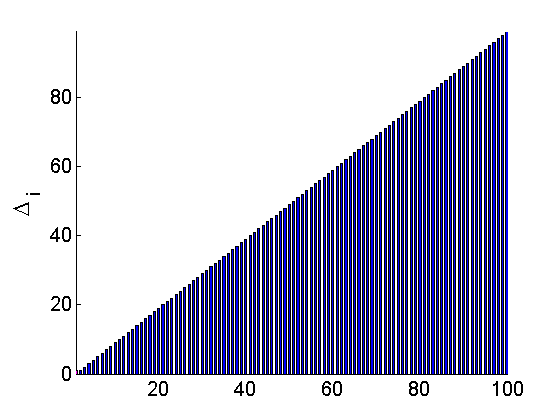}\qquad \includegraphics[width=0.33\textwidth,clip=true,trim=10 15 25 0]{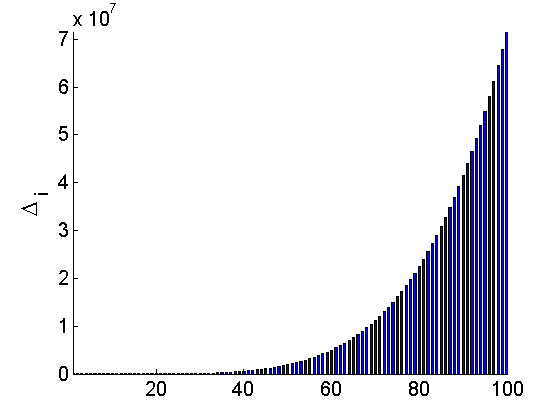}
  \caption{Constants in \eqref{eq.loose_bound_m} for problem $2$ with history length $m=1$ (left plot) and $m=5$ (right plot).}
  \label{fig.p2delta}
\efigure

\bfigure[H]
  \centering
  \includegraphics[width=0.47\textwidth,clip=true,trim=15 5 55 15]{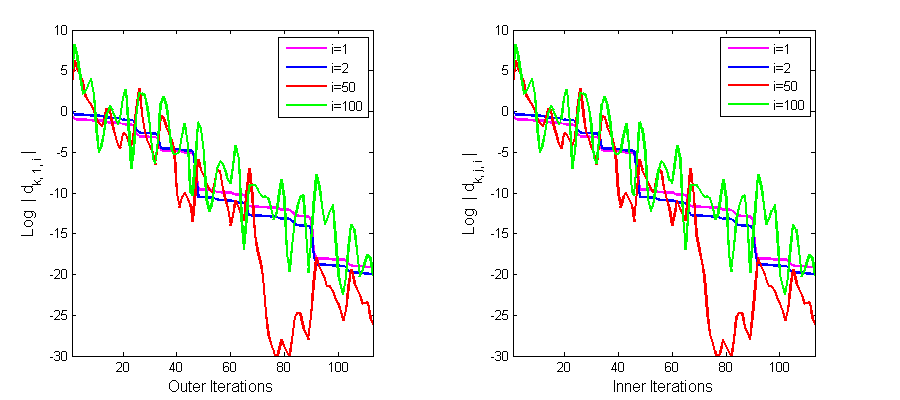}\qquad \includegraphics[width=0.47\textwidth,clip=true,trim=15 5 55 15]{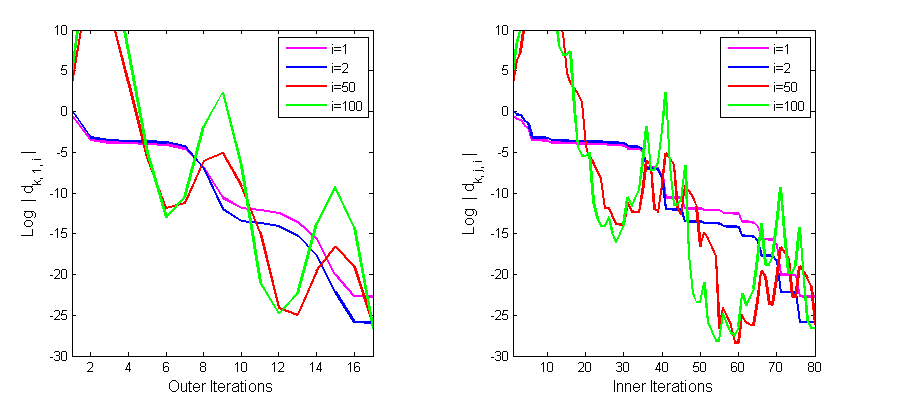}
  \caption{Weights in \eqref{eq.combination} for problem $3$ with history length $m=1$ (left two plots) and $m=5$ (right two plots).}
  \label{fig.p3}
\efigure

\bfigure[H]
  \centering
  \includegraphics[width=0.33\textwidth,clip=true,trim=10 15 25 0]{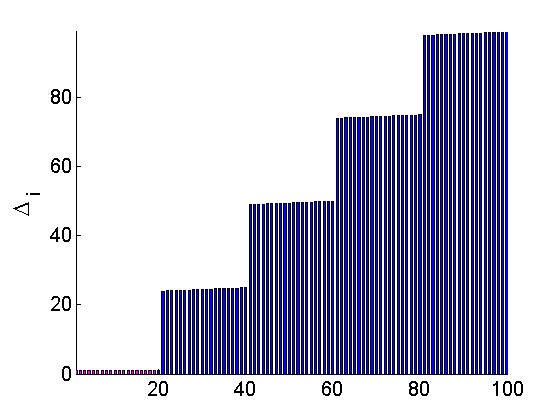}\qquad \includegraphics[width=0.33\textwidth,clip=true,trim=10 15 25 0]{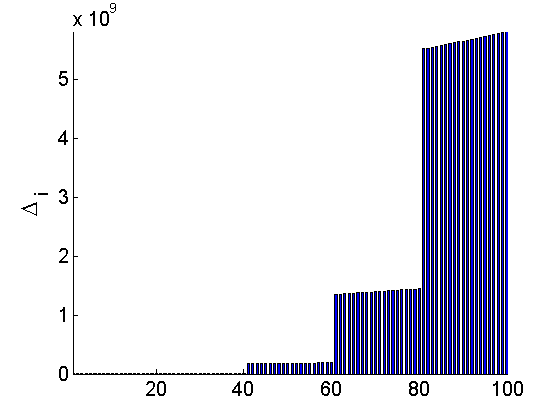}
  \caption{Constants in \eqref{eq.loose_bound_m} for problem $3$ with history length $m=1$ (left plot) and $m=5$ (right plot).}
  \label{fig.p3delta}
\efigure

\bfigure[H]
  \centering
  \includegraphics[width=0.47\textwidth,clip=true,trim=15 5 55 15]{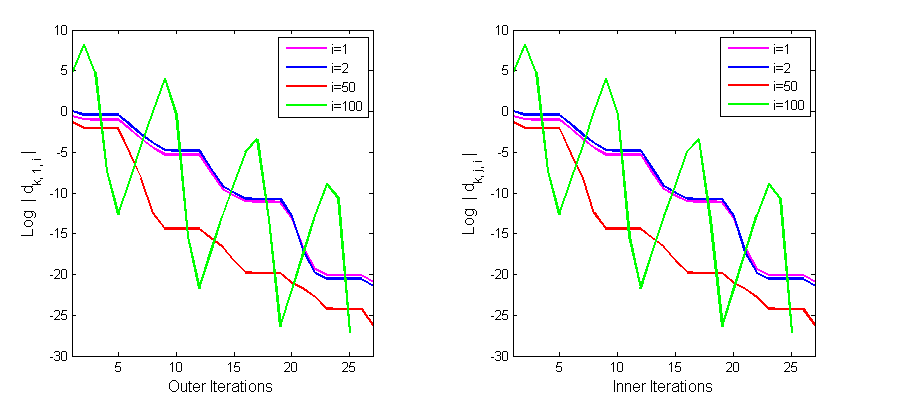}\qquad \includegraphics[width=0.47\textwidth,clip=true,trim=15 5 55 15]{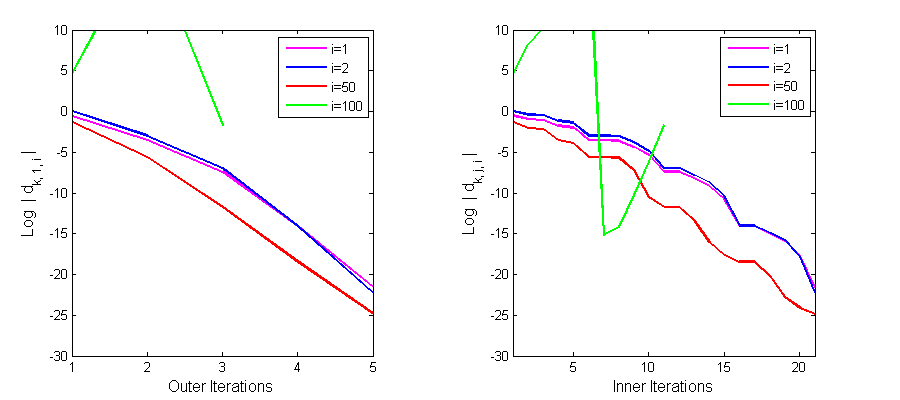}
  \caption{Weights in \eqref{eq.combination} for problem $4$ with history length $m=1$ (left two plots) and $m=5$ (right two plots).}
  \label{fig.p4}
\efigure

\bfigure[H]
  \centering
  \includegraphics[width=0.33\textwidth,clip=true,trim=10 15 25 0]{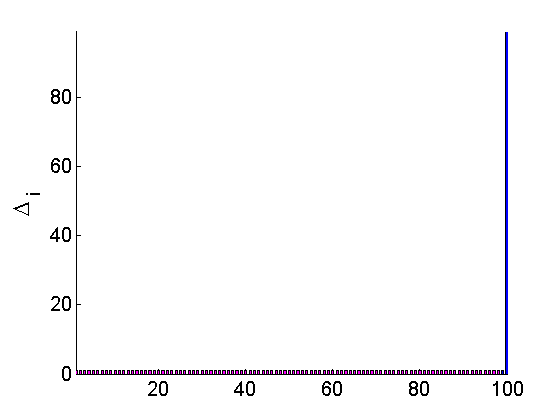}\qquad \includegraphics[width=0.33\textwidth,clip=true,trim=10 15 25 0]{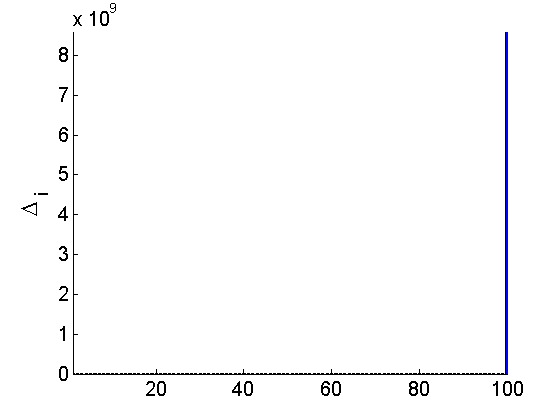}
  \caption{Constants in \eqref{eq.loose_bound_m} for problem $4$ with history length $m=1$ (left plot) and $m=5$ (right plot).}
  \label{fig.p4delta}
\efigure

\bfigure[H]
  \centering
  \includegraphics[width=0.47\textwidth,clip=true,trim=15 5 55 15]{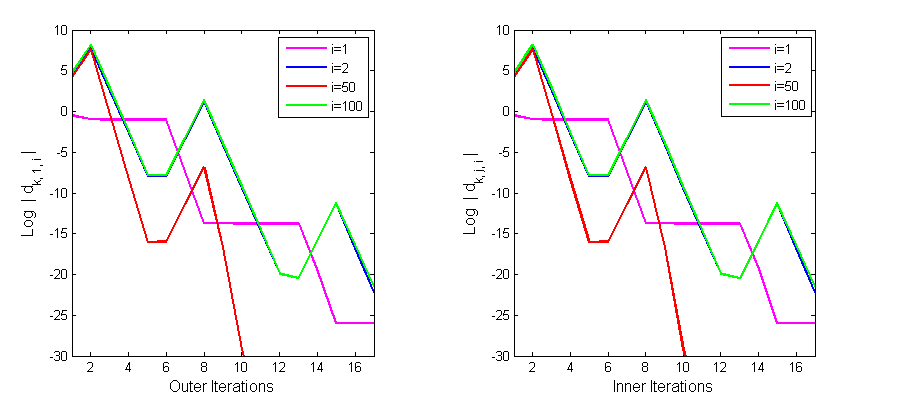}\qquad \includegraphics[width=0.47\textwidth,clip=true,trim=15 5 55 15]{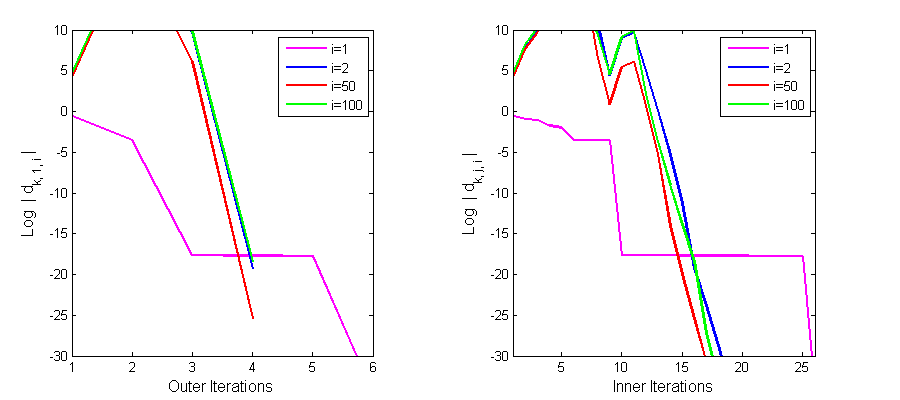}
  \caption{Weights in \eqref{eq.combination} for problem $5$ with history length $m=1$ (left two plots) and $m=5$ (right two plots).}
  \label{fig.p5}
\efigure

\bfigure[H]
  \centering
  \includegraphics[width=0.33\textwidth,clip=true,trim=10 15 25 0]{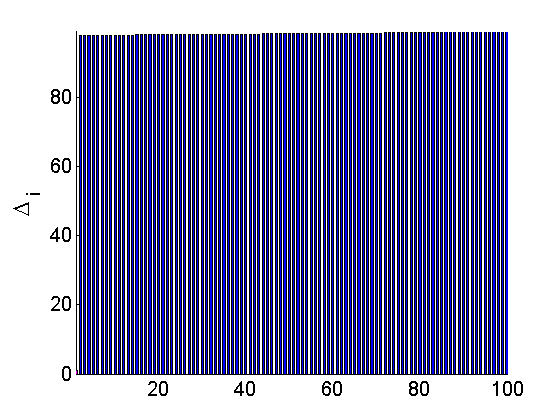}\qquad \includegraphics[width=0.33\textwidth,clip=true,trim=10 15 25 0]{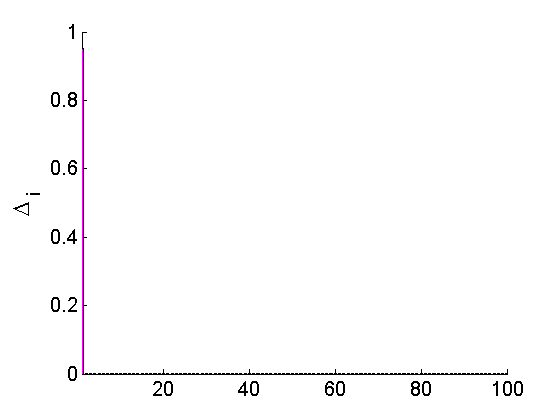}
  \caption{Constants in \eqref{eq.loose_bound_m} for problem $5$ with history length $m=1$ (left plot) and $m=5$ (right plot).}
  \label{fig.p5delta}
\efigure

\section{Conclusion}\label{sec.conclusion}

We have shown that the limited memory steepest descent (LMSD) method proposed by \cite{Flet12} possesses an $R$-linear rate of convergence for any history length $m \in \N{}$ when it is employed to minimize a strongly convex quadratic function.  Our analysis effectively extends that in \cite{DaiLiao02}, which covers only the $m=1$ case.  We have also provided the results of numerical experiments to demonstrate the practical performance of the algorithm, the results of which are informed by our theoretical analysis.


\ifthenelse{\coralreport = 1}{
\bibliographystyle{plain}
}{
\bibliographystyle{IMANUM-BIB}
}
\bibliography{lmsd_convergence}

\end{document}